\documentclass[preprint]{imsart}

\usepackage{mathtools}
\usepackage{paralist}

\usepackage{amsmath,amssymb,amsthm}
\RequirePackage[colorlinks,citecolor=cyan,urlcolor=blue,linkcolor=blue]{hyperref}
\usepackage[a4paper]{geometry}                % See geometry.pdf to learn the layout options. There are lots.
\usepackage{subfig,float}

\usepackage{enumerate}
\usepackage{graphicx}
\usepackage{verbatim,epsfig,color,lscape}
\usepackage{epstopdf}
\usepackage{color}
\usepackage{natbib}
\usepackage{algorithm}
\usepackage{algpseudocode}

\def\l{\left}
\def\r{\right}
\newcommand{\be}[1]{\begin{equation*}#1\end{equation*}}
\newcommand{\ben}[1]{\begin{equation}#1\end{equation}}

\def\ml#1{\begin{multline*}{#1}\end{multline*}}
\def\mln#1{\begin{multline}{#1}\end{multline}}

\newcommand{\wt}{\widetilde}
\newcommand{\wh}{\widehat}

\renewcommand{\exp}[1]{\operatorname{exp}\left(#1\right)} % Exponential
\newcommand{\h}[1]{h^{(#1)}}
\newcommand{\m}{\mathcal}
\newcommand{\mb}{\mathbb}
\newcommand\argmin{\mathop{\mbox{argmin}}}

\newcommand{\var}{\mbox{Var}}

\newcommand{\eps}{\varepsilon}

\newcommand{\med}[1]{\mbox{med}\left(#1\right)}

\newcommand{\pr}[1]{\mathbb{P}{\left(#1\right)}}

\newtheorem{lemma}{Lemma}
\newtheorem{theorem}{Theorem}
\newtheorem{remark}{Remark}

\begin{document}

\begin{frontmatter}
\title{Efficient median of means estimator}
\runtitle{Efficient MOM}

\begin{aug}
\author{\fnms{Stanislav} \snm{Minsker}\ead[label=e1,mark]{minsker@usc.edu} \thanksref{e1,a}}
\address[a]{Department of Mathematics, University of Southern California \\
\printead{e1}}
\thankstext{t3}{Author acknowledges support by the National Science Foundation grants DMS CAREER-2045068 and CCF-1908905.}
\end{aug}

\begin{abstract}%
The goal of this note is to present a modification of the popular median of means estimator that achieves sub-Gaussian deviation bounds with nearly optimal constants under minimal assumptions on the underlying distribution. We build on a recent work on the topic by the author, and prove that desired guarantees can be attained under weaker requirements.
\end{abstract}

%\maketitle

\end{frontmatter}

\mathtoolsset{showonlyrefs}
%################
\section{Introduction.}
\label{sec:intro}
%################

Let $X$ be a random variable with mean $\mu$ and variance $\sigma^2$. A sub-Gaussian estimator of $\mu$ based on a sample $\m X=\{X_1,\ldots,X_N\}$ of i.i.d. copies of $X$ is a measurable function $\wt\mu:=\wt \mu(\m X; t)$ such that 
$\pr{\l| \wt\mu -\mu \r|\geq C\sigma\sqrt{\frac t N}} \leq c e^{-t}$ for a absolute constants $c,C>0$ and all $t\in \l[1,t_{\mathrm{max}}(N)\r]$. It is known (for instance, see the work by \citet{catoni2012challenging}) that $C\geq \sqrt{2}$. A natural question, posed previously by \citet{devroye2016sub}, is whether sub-Gaussian estimators with $C=\sqrt{2}+o(1)$, where $o(1)$ is a function that goes to $0$ as $N$ (and possibly $t$) tend to infinity, exist. 

Several authors showed that such estimators can indeed be constructed under various additional assumptions. In one of the earliest works on the topic, \citet{catoni2012challenging} presented the first known example of sharp sub-Gaussian estimators for distributions with finite fourth moment and a known upper bound on the kurtosis, as well as for distributions with known variance. Construction by \citet{devroye2016sub} similarly required the fourth moment to be finite.
%however, it did not require the value of the kurtosis to be known.
% \citet{minsker2021robust} designed an asymptotically efficient sub-Gaussian estimator $\wt \mu_N$ that satisfies $\sqrt{N}\l( \wt\mu_N - \mu\r)\xrightarrow{d} N(0,\sigma^2)$ assuming only the finite second moment plus a mild, ``small-ball'' type condition. However, the constants in the non-asymptotic version of their bounds were not sharp. 
One of the strongest results is the one by \citet{lee2020optimal}: their estimator attains required guarantees uniformly over the class of distributions with finite variance, assuming just the finite second moment, albeit with $C=\sqrt{2}$ only in the limit as $t \to \infty$. \citet{minsker2022u} proposed a permutation-invariant version of the well known median of means (MOM) estimator \citep[][]{Nemirovski1983Problem-complex00,jerrum1986random,alon1996space} and proved that it achieves desired guarantees for the class of distributions with more than $\frac{3+\sqrt 5}{2}$ finite moments and ``sufficiently regular'' probability density functions.  

The main goal of this essay is to present a modification of the ``permutation-invariant'' MOM estimator that attains sub-Gaussian guarantees with asymptotically optimal constants for distributions possessing $2+\eps$ moments for some $\eps>0$. 
This result could yield improvements for a variety of robust algorithms (e.g., see the survey by \citet[][]{lugosi2019mean}) that rely on the classical MOM estimator serves as a subroutine. 

%##############################
\subsection{Notation.}
%##############################

For a positive integer $N$, $[N]$ will denote the set $\{1,\ldots,N\}$. 
We employ standard big-O and small-o notation for asymptotic relations between functions and sequences; it will be implicitly assumed that $o(1)$ and $O(1)$ may denote different functions from line to line. 
Moreover, given two sequences $\{a_n\}_{n\geq 1}$ and $\{b_n\}_{n\geq 1}$ where $b_n\ne 0$ for all $n$, we will write that $a_n\ll b_n$ if $\frac{a_n}{b_n}=o(1)$ as $n\to\infty$. Additional notation will be introduced in the main text whenever necessary.

%#################
\section{Main results.}
%#################

Let us recall the definition of the classical median of means estimator. 
Given an i.i.d. sample $\m X=\{X_1,\ldots,X_N\}$ from distribution $P$ with mean $\mu$ and variance $\sigma^2$, let $G_1\cup\ldots\cup G_k\subseteq [N]$ be a collection of $k$ disjoint subsets of cardinality $\lfloor N/k\rfloor$ each, $\bar X_j:=\frac{1}{|G_j|}\sum_{i\in G_j} X_i$ and 
$\wh\mu_{\mathrm{MOM}} = \med{\bar X_1,\ldots,\bar X_k},$ where $\med{\cdot}$ stands for the ``median.''
It is known that $\wh\mu_{\mathrm{MOM}}$ satisfies the inequality $\pr{\l| \wh\mu_{\mathrm{MOM}} -\mu \r|\geq C\sigma\sqrt{\frac t N}} \leq 2 e^{-t}$ with $C=\sqrt{\pi}+o(1)$, where $o(1)$ goes to $0$ as $k,N/k\to\infty$. 
\citet{minsker2022u} proved that allowing the overlapping subsets of data improves the constant: given $J\subseteq [N]$ of cardinality $|J|=\lfloor N/k\rfloor$, set $\bar X_J:= \frac{1}{|J|}\sum_{j\in J}X_j$ and define 
$\wh\mu_{U} =\med{\bar X_J, \ |J|=\lfloor N/k\rfloor}$, where $\l\{\bar X_J, \ |J|=\lfloor N/k\rfloor \r\}$ denotes the set of sample averages computed over all possible subsets of $[N]$ of cardinality $\lfloor N/k\rfloor$. Then $\wh\mu_U$ attains sub-Gaussian deviations with $C=\sqrt 2+o(1)$ under the assumptions described in section \ref{sec:intro}. Essentially, $\wh\mu_U$ is a function of the order statistics which are complete and sufficient for the family of all distributions with finite variance. 

Our construction, presented below, shows that it is not necessary to use all possible sample means, and that a much smaller collection of averages suffices: not only this makes computation easier, but the theoretical guarantees for the resulting estimator hold under weaker assumptions. The main idea is to split the data into subsets of size smaller than $\lfloor N/k\rfloor$, and construct all possible sample means using these subsets as ``building blocks''. The size of the overlap is then naturally proportional to the size of the block. For instance, the estimator $\wh\mu_U$ corresponds to the blocks of size $1$, resulting in the sample means over all possible subsets of a given size. Our results show that allowing the block size to be slowly growing with the the sample size could be beneficial. Formally, let $k,l$ be positive integers such that $\frac{\lfloor N/k\rfloor}{l} \in \mb N$. 
Assume that $G_1\cup\ldots\cup G_{lk}\subseteq [N]$ are disjoint subsets of cardinality $\lfloor \frac{N}{lk}\rfloor$ each, and $Z_j := \bar X_j =\frac{1}{|G_j|}\sum_{i\in G_j} X_i, \ j=1,\ldots,lk$. It will be convenient to set $n=lk$, $m=\lfloor \frac{N}{k}\rfloor$, and to view $Z_1,\ldots,Z_{n}$ is a new i.i.d. sample; clearly, $Z_1$ has mean $\mu$ and variance $\frac{\sigma^2}{m/l}$. 
Given $J\subseteq [n]$ of cardinality $|J|=l$, set $\bar Z_J:= \frac{1}{l}\sum_{j\in J}Z_j$; note that $\bar Z_J$ is an average of $m$ observations from the original sample $\m X_N$, same as in the definition of the standard MOM estimator. Define $\m A_{n}^{(l)} = \l\{ J\subset [n]: \ |J|=l\r\}$ and
\ben{
\label{eq:U-mom-est}
\wh\mu_N: =\med{\bar Z_J, \ J\in \m A_{n}^{(l)}},
}
where $\l\{\bar X_J, \ J\in \m A_{n}^{(l)}\r\}$ denotes the set of sample averages computed over all possible subsets of $[n]$ of cardinality $l$. In other words, $\wh\mu_N$ is the median of means computed over overlapping subsets of data, where the size of the overlap is proportional to $\lfloor N/lk\rfloor$, the size of the block $G_1$. 
We remark here that all explicit, non-asymptotic deviations guarantees that are valid for the classical MOM estimator $\wh\mu_{\mathrm{MOM}}$ automatically extend to $\wh\mu_N$ in view of the so-called ``Hoeffding representation'' of U-statistics \citep{lee1990u} as the average of averages of independent random variables; pursuit of optimal constant however appears to require the bounds that include asymptotic terms. Everywhere below, it is assumed that $k,m,l$ and functions of the sample size $N$. We proceed with the statement of our main result. Denote 
%necessitates reliance on
\ben{
\label{eq:g}
g(m):= \frac{6}{\sqrt m} \mb E\l[ \l(\frac{X_1 - \mu}{\sigma}\r)^2 \min\l(\l|\frac{X_1-\mu}{\sigma}\r|, \sqrt m\r) \r].
}
\citet{feller1968berry} proved that $g(m)$ controls the rate of convergence in the central limit theorem, namely that $\sup_{t\in \mb R}\l| \Phi_m(t) - \Phi(t)\r|\leq g(m)$ where $\Phi_m$ and $\Phi$ are the distribution functions of $\frac{\sum_{j=1}^m X_j - \mu}{\sigma\sqrt{m}}$ and the standard normal law respectively. It is well known that $g(m)\to 0$ as $m\to\infty$ for distributions with finite variance. Moreover, $g(m)$ admits an upper bound of the form $g(m) \leq C \mb E\l| \frac{X_1-\mu}{\sigma}\r|^{2+\eps} m^{-\eps/2}$ whenever $\mb E|X_1-\mu|^{2+\eps}<\infty$ for some $\eps\in(0,1]$. In the context of the median of mean estimation, the role of $g(m)$ is to control the difference between the mean and the median corresponding to the distribution of $\frac{1}{m}\sum_{j=1}^m X_j$, which can be seen as the main contribution to the the bias of $\wh\mu_N$.
\begin{theorem}
\label{th:main-1}
Assume that $\mb E|X_1-\mu|^{2+\eps}<\infty$ for some $\eps>0$. Suppose that $l=o(m^\eps)$ and let $L(n,l)$ and $M(n,l)$ be any sequences such that $L(n,l)\gg \frac{n}{l}g^2(m)$ and $M(n,l)\ll \frac{n}{l^2}$. 
% a function that is increasing in $n$, decreasing in $l$, and such that $q(n,l) = o(n/l^2)$. 
 Then for all $L(n,l)\leq t\leq M(n,l)$, 
\[
\pr{\l| \wh\mu_N - \mu\r|\geq \sigma\sqrt{\frac tN}}\leq 3\exp{-\frac{t}{2(1+o(1))}},
\]
where $o(1)\to 0$ as $l,k\to\infty$ uniformly over all $t\in \l[ L(n,l),M(n,l) \r]$.
\end{theorem}
\begin{remark} \hfill 
\begin{enumerate}
\item[(a)]
%(N/k)^\eps}
A possible choice of parameters is $l=\log(m)$, $L(n,l) = \frac{n}{l}\frac{\log(m)}{m^\eps}$ and $M(n,l) = \frac{n}{l^2\log(l)}$. By varying $k$, the deviation guarantees can be attained in the desired range of the confidence parameter.
\item[(b)] The question of uniformity of the bounds with respect to the underlying distribution is not explicitly addressed in this note. In particular, the $o(1)$ quantities appearing in the inequalities are distribution-dependent. With additional effort, it should be possible to prove uniformity with respect to the classes of distributions $\m P_N$ of $X$ satisfying moment conditions of the form $\mb E \l|\frac{X-\mu}{\sigma}\r|^{2+\eps}\leq a_N$ for a sequence $a_N$ that grows sufficiently slow.
\item[(c)] Exact computation of $\wh\mu_N$ is still prohibitively expensive from a numerical standpoint, as the naive upper bound for evaluating the estimator exactly is $O\l( (n/l)^l \log (n/l)\r)$. Instead, one may select a collection of $T$ subsets among $J\in \m A_n^{(l)}$ uniformly at random and compute the median of the corresponding sample means: in view of Theorem 1 in section 4.3.3 of the book by \citep{lee1990u} implies that the \emph{asymptotic} distribution of the estimator constructed in this way coincides with the asymptotic distribution $N(0,\sigma^2)$ of $\wh\mu_N$ as soon as $T\gg n/l$. However, this asymptotic equivalence does not automatically imply sharp non-asymptotic bounds of the estimator computed from subsampled blocks any more: results of such nature are currently unknown to us and require further investigation.
\end{enumerate}
\end{remark}

\begin{proof}
As $\wh\mu_N$ is scale-invariant, we can and will assume that $\sigma^2=1$. 
Set $\rho(x) = |x|$, and note that the equivalent characterization of $\wh\mu$ as an M-estimator is
\[
\wh\mu \in \argmin_{z\in \mb R} \sum_{J\in \m A_{n}^{(l)}} \rho\l( \sqrt{m}\l( \bar Z_J - z\r)\r).
\]
The necessary conditions for the minimum of $F(z) := \sum_{J\in \m A_{n}^{(l)}}  \rho\l( \sqrt{m}\l( \bar Z_J - z\r)\r)$ imply that $0\in \partial F(\wh\mu_N)$, hence the left derivative 
$F'_-(\wh\mu_N)\leq 0$. Therefore, if $\sqrt{N}\l(\wh\mu_N - \mu\r)\geq \sqrt t$ for some $t > 0$, then $\wh\mu_N \geq \mu + \sqrt{t/N}$ and, due to $F'_-$ being nondecreasing, $F'_-\l( \mu + \sqrt{t/N} \r)\leq 0$. It implies that 
\mln{
\label{eq:b001}
\pr{\sqrt{N}(\wh\mu_N - \mu)\geq \sqrt t} \leq \pr{ \sum_{J\in \m A_{n}^{(l)}}   \rho'_-\l( \sqrt{m}\l(  \bar Z_J - \mu - \sqrt{t/N}\r)\r) \geq 0} 
\\
= \pr{\frac{\sqrt{k}}{{n\choose l}} \sum_{J\in \m A_{n}^{(l)}} \l(\rho'_-\l( \sqrt{m}\l( \bar Z_J - \mu - \sqrt{t/N}\r)\r) - \mb E\rho'_-\r) \geq 
-\sqrt{k}\mb E\rho'_- }
}
where we used the shortcut $\mb E\rho'_-$ in place of 
%\begin{align*}
\ml{
\mb E\rho'_-\l( \sqrt{m}\l( \bar Z_J - \mu - \sqrt{t/N}\r)\r) = I\l\{ \sqrt{m}(\bar Z_J - \mu) \leq \sqrt{\frac{t}{k}}\r\} - I\l\{ \sqrt{m}(\bar Z_J - \mu) > \sqrt{\frac{t}{k}}\r\}
\\
= 1 - 2 I\l\{ \sqrt{m}(\bar Z_J - \mu) \leq \sqrt{\frac{t}{k}}\r\}.
}
%\end{align*}
Note that 
\mln{
\label{eq:b11}
-\sqrt{k} \mb E\rho'_-\l( \sqrt{m}\l(\bar Z_J - \mu - \sqrt{t/N}\r)\r) = 
-\sqrt{k}\l(1 - 2 \pr{ \sqrt{m}\l( \bar Z_J - \mu - \sqrt{t/N}\r) \leq 0} \r)
\\
= 2\sqrt{k}\l( \Phi\l( \sqrt{\frac{ t}{k}}\r) - \Phi(0) \r) - 2\sqrt{k}\l(\Phi\l( \sqrt{\frac{t}{k}} \r) - \pr{ \sqrt{m}\l( \bar Z_J - \mu \r) \leq \sqrt{\frac{t}{k}}} \r)
\\
\geq -2\sqrt{k} \cdot g(m) + 2\sqrt t \frac{1}{\sqrt t/\sqrt{k}}\l( \Phi\l( \frac{\sqrt t}{\sqrt{k}}\r) - \Phi(0) \r). 
}
Since 
\ben{
\label{eq:b12}
2 \sqrt{t}\frac{1}{\sqrt t/\sqrt{k}}\l( \Phi\l( \frac{\sqrt t}{\sqrt{k}}\r) - \Phi(0) \r) = 2 \sqrt t \l( \phi(0) + O(\sqrt{t/k}) \r)
= \sqrt t \l(\sqrt{\frac 2 \pi} + O(\sqrt{t/k}) \r)
}
where $\phi(t) = \Phi'(t)$, we see that 
\[
-\sqrt{k}\, \mb E\rho'_-\l( \sqrt{m}\l( \bar Z_J - \mu - \sqrt{t/N}\r)\r) 
\geq -2\sqrt{k}\cdot g(m) + \sqrt t\l( \sqrt{\frac 2 \pi} + O(\sqrt{t/k}) \r)
\]
which is $\sqrt t\sqrt{\frac 2 \pi} \l(1 + o(1)\r)$ whenever $t\ll k$ and $t \gg k\, g^2(m)$. 
It remains to analyze the U-statistic 
\ben{
\label{eq:u-stat1}
\sqrt{k} \, U_{n,l}(\rho'_-) = \frac{\sqrt{k}}{{n\choose l}} \sum_{J\in \m A_{n}^{(l)}} \l(\rho'_-\l( \sqrt{m}\l( \bar Z_J - \mu - \sqrt{t/N}\r)\r) - \mb E\rho'_-\r).
}
As the expression above is invariant with respect to the shift $Z_j \mapsto Z_j - \mu$, we can assume that $\mu=0$. 
For $i \in [N]$, let 
\be{
h^{(1)}(Z_i) 
=\sqrt{l}\, \mb E\l[ \rho'_-\l(\sqrt m\l( \frac{1}{l}\sum_{j=1}^{l-1} \tilde Z_j + \frac{Z_i}{l} - \sqrt{t/N}\r)\r)\,\big|\,Z_i \r] - \sqrt{l}\,\mb E\rho'_-,
}
where $(\tilde Z_1,\ldots,\tilde Z_l)$ is an independent copy of $(Z_1,\ldots,Z_l)$ based on a sample $\tilde{\m X}_N$ that is an independent copy of $\m X_N$. Our goal is to determine the size of $\var(h^{(1)}(X_1))$.
\begin{remark}
The quantity $h^{(1)}(Z)$ is related to the so-called H\'{a}jek projection that can be viewed as the best (in mean squared sense) approximation of the U-statistic $U_{n,l}(\rho'_-)$ in terms of the sums of i.i.d. random variables. For related background on U-statistics, we refer the reader to an excellent monograph by \citet{lee1990u}.
\end{remark}

\begin{lemma}
\label{lemma:var1}
In the framework of Theorem \ref{th:main-1}, 
\[
\var\l(h^{(1)}(Z_1) \r) \to \frac{2}{\pi}
\] 
as $l,k\to \infty$, uniformly over all $t\in \l[ L(n,l),M(n,l) \r]$.
\end{lemma}
The proof of the lemma is given in section \ref{proof:A}. 
The following result, a deviation inequality for U-statistics of order that grows with the sample size, is the second key technical tool required to complete the argument. 

\begin{theorem}
\label{th:concentration}
Let $h:\mb R^l\mapsto \mb R$ be a function that is invariant with respect to permutations of its arguments, and let 
$U_{n,l}(h) = \frac{1}{{n\choose l}} \sum_{J\in \m A_{n}^{(l)}} \l(h\l(X_j, \ j\in J\r)  - \mb Eh\l(X_1,\ldots,X_l\r)\r)$ be the corresponding U-statistic with kernel $h$ evaluated on a sample $X_1,\ldots,X_n$. Assume that $l$ is an increasing function of $n$, and that 
\begin{enumerate}
\item[(a)] $h$ is uniformly bounded;
\item[(b)] $\liminf_{l \to\infty}\var\l( \sqrt{l}\, h^{(1)}(X_1) \r)>0$, where $h^{(1)}(X_1) = \mb E\l[ h(X_1,X_2,\ldots,X_l)\vert X_1\r]$.    
\end{enumerate}
Let $q(n,l)$ be increasing in $n$, decreasing in $l$, and such that $q(n,l) = o\l( \frac{n}{l^2}\r)$. 
Then for all $2\leq t\leq q(n,l)$,
\[
\pr{\l|U_{n,l}(h)\r| \geq \sqrt{\frac{tl}{n}}} \leq (2+o(1))\exp{-\frac{t}{2(1+o(1))\var\l( \sqrt{l}\, h^{(1)}(X_1) \r)}},
\] 
where $o(1)\to 0$ as $l, n/l\to \infty$ uniformly over $2\leq t\leq q(n,l)$. 
\end{theorem}
The proof of this result is outlined in section \ref{proof:B} \footnote{We note that closely related results for U-statistics were obtained by \citet{maurer2019bernstein}, and it may be possible to use Maurer's inequality in place of Theorem \ref{th:concentration}.}. To get the desired inequality for the estimator $\wh\mu_N$, it remains to apply Theorem \ref{th:concentration} and Lemma \ref{lemma:var1} to the U-statistic defined in \eqref{eq:u-stat1}: specifically, we deduce that 
\[
\pr{\l|\sqrt{k} \, U_{n,l}(\rho'_-)\r| \geq \sqrt t\sqrt{\frac 2 \pi} \l(1 + o(1)\r) }\leq 2\exp{-\frac{t}{2(1+o(1))}},
\]
uniformly over $\frac{n}{l}g^2(m)\ll t \ll \frac{n}{l^2}$, and the final result follows.

\end{proof}

%################################
\section{Proof of Lemma \ref{lemma:var1}.}
\label{proof:A}
%################################

Note that we can rewrite $h^{(1)}(Z_1)$ as 
\[
h^{(1)}(Z_1) 
=\sqrt{l}\, \mb E\l[ \rho'_-\l(\sqrt m\l( \frac{1}{m}\sum_{j=1}^{m-m/l} \tilde Z_j  + \frac{1}{\sqrt{ml}}\l( Z_1\sqrt{m/l}\r) - \sqrt{t/N}\r)\r)\,\big|\,Z_i \r] - \sqrt{l}\,\mb E\rho'_-.
\]
Given an integer $r\geq 1$, let $\wt\Phi_{r}(t)$ be the cumulative distribution function of $\sum_{j=1}^r \tilde X_j$. Then 
\ml{
h^{(1)}(Z_1) = \sqrt{l}\l(2\wt\Phi_{m-m/l}\l(m\sqrt{\frac{t}{N}} - \sqrt{m/l}\l( Z_1\sqrt{m/l}\r)\r) -1 \r) - \sqrt{l}\, \mb E\,\rho'_- 
\\
= 2\sqrt{l}\l( \wt\Phi_{m-m/l}\l( m\sqrt{\frac{t}{N}} - \sqrt{m/l}\l( Z_1\sqrt{m/l}\r)\r)  - 
\mb E \wt\Phi_{m-m/l}\l( m\sqrt{\frac{t}{N}} - \sqrt{m/l}\l( Z_1\sqrt{m/l}\r)\r) \r)
\\
=2\sqrt{l}\int_\mb R \l( \wt\Phi_{m-m/l}\l( m\sqrt{\frac{t}{N}} - \sqrt{m/l}\l( Z_1\sqrt{m/l}\r)\r) \r.
\\
\l. - \wt\Phi_{m-m/l}\l( m\sqrt{\frac{t}{N}} - x\sqrt{m/l} \r) \r) dP_{Z_1\sqrt{m/l}}(x),
}
with $P_{Z_1\sqrt{m/l}}$ being the law of $Z_1\sqrt{m/l}$. Feller's version of Berry-Esseen theorem implies that 
\[
\sup_{x\in \mb R}\l| \wt\Phi_{m-m/l}(x) - \Phi(x/\sqrt{m-m/l}) \r|\leq 6g(m-m/l)
\] 
where $\Phi$ is the distribution function of standard normal law. Therefore, 
\ml{
\l| h^{(1)}(Z_1) - 2\sqrt{l}\int_\mb R 
\l(\Phi\l( \frac{m}{\sqrt{m-m/l}}\sqrt{\frac{t}{N}} - \sqrt{\frac{m/l}{m-m/l}}\l( Z_1\sqrt{m/l}\r)\r) \r.\r.
\\
\l.\l. - \Phi\l( \frac{m}{\sqrt{m-m/l}}\sqrt{\frac{t}{N}} - x\sqrt{\frac{m/l}{m-m/l}} \r) \r) dP_{Z_1\sqrt{m/l}}(x) \r| 
\leq 12\sqrt{l} \,g(m-m/l) \to 0
}
by assumption. At the same time, 
\ml{
\sqrt{l}\l(\Phi\l( \frac{m}{\sqrt{m-m/l}}\sqrt{\frac{t}{N}} - \sqrt{\frac{m/l}{m-m/l}}\l( Z_1\sqrt{m/l}\r)\r) - \Phi\l( \frac{m}{\sqrt{m-m/l}}\sqrt{\frac{t}{N}} - x \sqrt{\frac{m/l}{m-m/l}} \r) \r)
\\ 
= \frac{1}{\sqrt{2\pi}}\exp{-q(x)/2} (x - Z_1\sqrt{m/l})\sqrt{\frac{m}{m-m/l}} + \frac{C(x,Z_1)}{\sqrt{l}} (x - Z_1\sqrt{m/l})^2 \frac{m/l}{m-m/l}
}
where $q(x):=\l( \frac{m}{\sqrt{m-m/l}}\sqrt{\frac{t}{N}} - x\sqrt{\frac{m/l}{m-m/l}}\r)^2$ is such that $q(x)\to 0$ as $l\to\infty$ and $C(x,Z_1)$ is a bounded function. Therefore, $h^{(1)}(Z_1) - \sqrt{\frac{2}{\pi}}Z_1\sqrt{m/l} \to 0$ almost surely, assuming that $\sqrt{l}g(m)=o(1)$. 
Finally, note that 
\ml{
\sqrt{l} \l( \wt\Phi_{m-m/l}\l( m\sqrt{\frac{t}{N}} - \sqrt{m/l}\l( Z_1\sqrt{m/l}\r)\r) - \wt\Phi_{m-m/l}\l( m\sqrt{\frac{t}{N}} - x\sqrt{m/l} \r) \r) 
\\
\leq \sup_{z}\sqrt l\,\pr{\sum_{j=1}^{m-m/l} \tilde X_j \in \Big(z, z+ \l| \sqrt{m/l}\l(x - Z_1\sqrt{m/l}\r)\r| \Big]} 
\leq C \l|x-Z_1\sqrt{m/l}\r|,
}
where the last inequality follows from the well known bound for the concentration function (Theorem 2.20 in the book by \cite{petrov1995limit}); here, $C=C(P)>0$ is a constant that may depend on the distribution of $X_1$. We therefore conclude that the sequence $\l( h^{(1)}(Z_1) - \sqrt{\frac{2}{\pi}}Z_1\sqrt{m/l}\r)^2$ is uniformly integrable (as $Z_1\sqrt{m/l}$ is), hence the claim follows. 

%########################
\section{Proof of Theorem \ref{th:concentration}.}
\label{proof:B}
%########################

The union bound together with Hoeffding's decomposition entails that for any $t>0$ and $0<\eps<1$ (to be chosen later as a decreasing sequence $\eps(l)$), 
\ml{
\pr{\l|U_{n,l}(h)\r| \geq \sqrt{\frac{tl}{n}}} 
\\
\leq \pr{\l| \frac{l}{n}\sum_{j=1}^n h^{(1)}(Z_j) \r|\geq (1-\eps) \sqrt t\sqrt{\frac{l}{n}}}
+ \pr{\l| \sum_{j=2}^l \frac{{l\choose j} }{{n\choose j}} \sum_{J\in \m A_n^{(j)}} \h{j}(Z_i,\, i\in J) \r| \geq \eps\sqrt t\sqrt{\frac{l}{n}}},
}
where $h^{(j)}, \ 2=1,\ldots,l$ are the degenerate kernels corresponding to the higher-order terms of Hoeffding's decomposition (not to be confused with the derivatives!). Specifically, 
\[
h^{(j)}(y_1,\ldots,y_j)=(\delta_{y_1} - P_Y)\times\ldots\times(\delta_{y_j}-P_Y)\times P_Y^{m-j} h,
\]
where $\delta_y$ is the point measure concentrated at $y$; in particular, $\delta_y(h) = h(y)$. It is known that $h^{(j)}$ can be viewed geometrically as orthogonal projections of $h$ onto a particular subspace of $L_2(P_Y^{m})$. 
We refer the reader to the book by \citet{lee1990u} for futher details related to the background material on U-statistics and the Hoeffding's decomposition. Bernstein's inequality yields that
\ml{
\pr{\l| \frac{l}{n}\sum_{j=1}^n h^{(1)}(Z_j) \r|\geq (1-\eps)\sqrt t\sqrt{\frac{l}{n}}} 
\\
\leq 2\exp{-\frac{(1-\eps)^2 \,t/2}{ \var\l( \sqrt{l}\, h^{(1)}(Z_1) \r) + (1-\eps)\frac13\sqrt{\frac{l}{n}}\|h\|_\infty t^{1/2}}} 
\\
=2\exp{-\frac{(1-\eps)^2\, t}{2\,\var\l( \sqrt{l}\, h^{(1)}(Z_1) \r)(1+o(1))}}
}
where $o(1)\to 0$ as $n/l\to\infty$ uniformly over $t$. 
It remains to control the expression involving higher order Hoeffding decomposition terms. To this end, we will show that it is bounded from above by $\exp{-\frac{t}{2\,\var\l( \sqrt{l}\, h^{(1)}(X_1) \r)}} \cdot o(1)$ where $o(1)\to 0$ uniformly over the range of $t$. To this end, we will need concentration inequality for the U-statistics of growing order established in \citet[][Theorem 4.1]{minsker2022u}. Set $t_{j,\eps} = \l( \eps\frac{\sqrt{t}}{j^2} \l( \frac{n}{l} \r)^{\frac{j-1}{2}}\r)^2$, and note that, in view of the union bound,
\ml{
\pr{\l| \sum_{j=2}^l \frac{{l\choose j} }{{n\choose j}} \sum_{J\in \m A_n^{(j)}} \h{j}(Z_i,\, i\in J) \r| \geq \eps\sqrt{t}\sqrt{\frac{l}{n}}}
\\ 
\leq \sum_{j=2}^l \pr{\l| \frac{{l\choose j} }{{n\choose j}} \sum_{J\in \m A_n^{(j)}} \h{j}(Z_i,\, i\in J) \r| \geq \sqrt{t_{j,\eps}}\sqrt{\frac{l}{n}}  }
\\
\leq l \max_{2\leq j\leq l} \exp{-c\min\l((t\eps^2)^{1/j}\l( \frac{n}{l}\r)^{\frac{j-1}{j}},\l( \frac{t\eps^2}{\|h\|_\infty^2}\r)^{\frac{1}{j+1}}  \l( \frac{nj}{l^2} \r)^{\frac{j}{j+1}} \r)},
}
where the last inequality follows from the first bound of Theorem 4.1 in \citet[][]{minsker2022u}. Whenever $l\log^2(l)\ll k$ and $\eps\gg \log^{-1/2}(l)$, the last expression is at most 
\[
\max_{2\leq j\leq l}\exp{-c_1\min\l((t\eps^2)^{1/j}\l( \frac{n}{l}\r)^{\frac{j-1}{j}},\l( \frac{t\eps^2}{\|h\|_\infty^2}\r)^{\frac{1}{j+1}}  \l( \frac{nj}{l^2} \r)^{\frac{j}{j+1}} \r)}.
\]
In turn, it is bounded by $e^{-\frac{c_2t}{\eps}}$ whenever $t<\frac{n}{l^2} \eps^4$. Desired conclusion follows. 

\section*{Acknowledgements.}
The author is grateful to the anonymous Referees for their insightful comments and constructive feedback that helped improve the quality of the presentation.
% Acknowledgments---Will not appear in anonymized version
%\acks{To be added.}

\bibliography{MoM,robustERM}
\bibliographystyle{apalike}
%\appendix

% \crefalias{section}{appendix} % uncomment if you are using cleveref

\end{document}